\documentclass[10pt]{siamonline1116}

\usepackage{slashbox}
\usepackage[utf8]{inputenc}
\usepackage{bm}
\usepackage{graphicx}
\usepackage{amsmath}
\usepackage[parfill]{parskip}
\usepackage{xcolor}
\usepackage[english]{babel}
\usepackage{amsfonts}
\usepackage{amssymb}
\usepackage{amsbsy}
\usepackage{hyperref}
\usepackage{epstopdf}
\ifpdf
  \DeclareGraphicsExtensions{.eps,.pdf,.png,.jpg}
\else
  \DeclareGraphicsExtensions{.eps}
\fi

\numberwithin{theorem}{section}

\newcommand{\D}[1]{{#1}'}
\newcommand{\B}[1]{\mbox{\boldmath${#1}$\unboldmath}}
\newcommand{\R}{\mathbb{R}}
\newcommand{\X}{\mathcal{X}}

\newcommand{\TheTitle}{Shadowing-based data assimilation method for partially observed models} 
\newcommand{\TheAuthors}{De Leeuw, Dubinkina}

\ifpdf
\hypersetup{
  pdftitle={\TheTitle},
  pdfauthor={\TheAuthors}
}
\fi

\headers{\TheTitle}{\TheAuthors}

\title{{\TheTitle}\thanks{Submitted to the editors \today.
\funding{The work of the first author was
partially supported by the research program Mathematics of Planet Earth 2014 EW project 657.014.001, which is
financed by the Netherlands Organisation for Scientific Research (NWO).}
}}

\author{
Bart de Leeuw\thanks{Centrum Wiskunde \& Informatica, PO Box 94079, 1090 GB Amsterdam, Netherlands, 
	(\email{b.m.de.leeuw@cwi.nl}, \email{s.dubinkina@cwi.nl}).} 
\and 
Svetlana Dubinkina\footnotemark[2]
}

\usepackage{amsopn}

\begin{document}
\maketitle

\begin{abstract}
In this article we develop further an algorithm for data assimilation based upon a shadowing refinement technique [de Leeuw et al., \emph{SIAM J. Appl. Dyn. Sys.}, 17 (2018), pp. 2446–2477] to take partial observations into account. 
Our method is based on regularized Gauss-Newton method.
We prove local convergence to the solution manifold and provide a lower bound on the algorithmic time step. 
We use numerical experiments with the Lorenz 63 and Lorenz 96 models to illustrate convergence of the algorithm
and show that the results compare favourably with a variational technique---weak-constraint four-dimensional variational method---and a shadowing technique--pseudo-orbit data assimilation. Numerical experiments show that a preconditioner chosen based on a cost function allows the algorithm to find an orbit of the dynamical system in the vicinity of the true solution.
\end{abstract}



\section{Introduction}
Data assimilation (DA) methods combine orbits from a dynamical model
with measurement data to obtain an improved estimate for the state of a physical system~\cite{Ja70}. 
Well known strong-constraint four-dimensional variational data assimilation (4DVar) aims at finding the optimal initial condition for the dynamical model 
such that the distance to observations is minimized under a constraint of the estimate being an orbit of the dynamical model~\cite{Ta97}. A drawback of strong-constraint 4DVar is that the number of local minima of the corresponding cost function increases dramatically with assimilation window---time window over which observations are assimilated into the dynamical model~\cite{Be91,MiGhGa94,PiVaTa96}.
An existing remedy in 4DVar is introduction of a model error term in the cost function and is called weak-constraint 4DVar (WC4Var)~\cite{Sa70, Tr06}.
Then an estimate is a pseudo-orbit of the dynamical model rather than an orbit.
An orbit satisfies PDE of a dynamical model exactly, while a pseudo-orbit up to a small $\varepsilon$. It has been shown in e.g.~\cite{Tr06} that
WC4Var allows longer assimilation windows compared to the strong-constraint 4DVar.

An alternative DA approach that allows long assimilation windows is based on a model having a shadowing property. Let $F$ be the exact time-$\Delta t$ flow map of an autonomous ODE $\dot{x}=f(x)$. 
Suppose $\B{u}$ is an $\varepsilon$-orbit in a neighborhood of a hyperbolic set for $F$.
If the components of $\B{u}$ are the iterates of a numerical integrator with local truncation error bounded by $\varepsilon$, 
then these define an $\varepsilon$-orbit of $F$. 
The shadowing lemma (e.g.~Theorem 18.1.2 of \cite{KaHa95}) states that, for every $\delta>0$ there exists $\varepsilon>0$ such that 
$\B{u}$ is $\delta$-shadowed by an orbit of $F$.
Shadowing refinement~\cite{GHYS90} employs the pseudo-orbit as an initial guess  and
iteratively refines the pseudo-orbit to obtain an improved approximation of a true solution. 
The inverse problem to shadowing is to determine an optimal initial condition $u_0$ for a numerical integration, 
such that the numerical iterates $\B{u}$ $\delta$-shadow a desired orbit of $\dot{x}=f(x)$.
\par

There exist several shadowing-type DA methods. A pseudo-orbit DA method (PDA)~\cite{DuSm14} and a noise reduction algorithm~\cite{BrPar01} seek a (pseudo-)trajectory of a dynamical model by minimizing a cost function. Local convergence to the solution manifold corresponding to $\dot{x} -f(x) = 0$ was proven for a class of iteration schemes assuming full observations~\cite{BrPar01}. In numerical experiments, the noise reduction algorithm uses the Laplace operator and PDA---an algorithmic time step to achieve convergence to
the solution manifold, though without a robust answer whether these are the good choices for the convergence.  
 Another shadowing-type DA method instead of minimizing a cost function, seeks zeros of a cost operator~\cite{Leetal18}. Obviously, the (nonunique) global minimum of the cost function is zero and this value is reached if and only if the corresponding cost operator is zero. 
 
A shortcoming of existing shadowing-type DA methods is that for initialization they use \emph{full} observations in space (or partial observations combined with an estimation obtained from another DA method). Up to now truly partial observations (without any preprocessing involving another DA method) have not been thoroughly considered in shadowing-type DA methods. 
Therefore, in this paper we consider an initial guess for a shadowing-type DA method that consists of partial observations 
and a background trajectory, which was obtained from model propagation starting at an arbitrary initial condition and \emph{without} DA.
We develop further the shadowing-based DA method~\cite{Leetal18} to account for partial observations
based on Levenberg-Marquardt regularization~\cite{Le44,Ma63}, and prove local convergence following~\cite{BrPar01}.
The Levenberg-Marquardt algorithm can be seen as a regularization of the Gauss-Newton method, which is used in the shadowing-based DA method of~\cite{Leetal18}. A regularization parameter controls algorithmic time step, making the Gauss-Newton method convergent to the solution manifold independently of the starting point. 
The Levenberg-Marquardt regularization is well used in nonlinear optimization and data assimilation in particular, 
e.g. variational data assimilation~\cite{Maetal16}, and ensemble Kalman filter~\cite{ChOl13}.\par

Despite being convergent to the solution manifold, a shadowing-type DA method might poorly approximate the true solution due to observations being used only as initial guess. Therefore, in this paper we introduce a preconditioner for the corresponding gradient flow that modifies the direction of the search such that the estimate remains in the vicinity of observations. This is done in the spirit of trust region methods~\cite{NoWr06}, which together with Gauss-Newton type methods, have been an inspiration for new algorithms to solve nonlinear least-squares problems, see e.g.~\cite{StKi11}.

The rest of the paper is organized as follows. 
 In Section~\ref{sec:NR}, we briefly recall the shadowing-based DA method for full observations. 
 In Section~\ref{sec:SDA}, we introduce the shadowing-based DA method for partial observations and prove local convergence.  
 In Section~\ref{sec:numerics}, we present results for the Lorenz 63 and the Lorenz 96 models. Finally, we draw the conclusions in Section~\ref{sec:con}.

\section{Noise reduction}\label{sec:NR}
We consider a discrete deterministic model 
\begin{equation}\label{eq:map}
	x_{n+1} = F_n (x_n), \quad x_n \in \R^m, \quad n=0,\dots,N-1,
\end{equation}
where $F_n: \R^m \to \R^m$. We assume $F_n$ to be $\mathcal{C}^3$ for all $n$. 
In many applications the model is defined by the time-discretization of an ordinary differential equation
$\dot{x} = f (t,x)$, $x(t) \in \R^m$, 
which in turn may be defined as the space-discretization of a partial differential equation (or system of PDEs).

Let the sequence $\B{X}:=\{\X_0,\dots,\X_{N}\}$ be a distinguished orbit of~\eqref{eq:map}, 
referred to as the true solution of the model, and presumed to be unknown.
Suppose we are given a sequence of partial noisy observations $\B{y}:=\{y_0,\dots y_{N}\}$ related to $\B{X}$ via 
\begin{equation*}
	y_n = H_n \X_n + \xi_n, \qquad y_n \in \mathbb{R}^d, \quad n=0,\dots,N,
\end{equation*}
where $H_n:\R^m \rightarrow \R^d$, $d\le m$, is the linear observation operator, and the noise variables $\xi_n$ are drawn from 
a normal distribution  $\mathcal{N}(0,R_n)$ with zero mean and known observational error covariance matrix $R_n$.

Data assimilation is the problem of finding a pseudo-orbit $\B{u} = \{ u_0, u_1, \dots, u_N\}$, $u_n \in \R^d$, of the model \eqref{eq:map}, such that the differences $\| y_n - H u_n\|$ and  $\| u_{n} - F_n(u_{n-1})\|$, $n=1,\dots,N$ are small in an appropriately defined sense. This is done with the aim of minimizing the unknown error $\| u_n - \X_n \|$; see for example \cite{Ta97,LaStZy15}. 
Well known WC4DVar aims at finding the optimal initial condition $u_0$ of \eqref{eq:map} to minimize a cost function
\[
 C_\mathrm{var}(u_0; \{ y_n \}) = 
\sum_{n=1}^N (y_n-Hu_n)^T R^{-1}(y_n-Hu_n) + (u_n - F_n(u_{n-1}))^TQ^{-1}(u_n - F_n(u_{n-1})),
\]
where the $Q$ is model error (see e.g.\cite{Sa70,LeDe85,TaCo87,Ta97} and references therein). 

Instead of minimizing a cost function, the shadowing-based DA method~\cite{Leetal18} searches for a zero of the \emph{cost operator}
\begin{equation}\label{eq:residual}
	G(\B{u}) = \begin{pmatrix} G_0(\B{u}) \\ G_1(\B{u}) \\ \vdots \\ G_{N-1}(\B{u}) \end{pmatrix},
	\qquad G_n(\B{u}) = u_{n+1} - F_n(u_n), \quad n=0,\dots,N-1,
\end{equation}
using a contractive iteration started from (a proxy of) complete, noisy observations. Therefore we call this method noise reduction DA method. This approach is motivated by research on numerical shadowing methods. We stress that, just as with strong-constraint 4DVar, noise reduction DA attempts to find an exact orbit of \eqref{eq:map} consistent with the observations.  However, instead of solving directly for the initial condition, we solve for the whole orbit at once.

Noise reduction DA seeks an update $\B{P}^{(k)}$ by approximately solving 
\begin{equation}\label{eq:Newton}
	G \left(\B{u}^{(k)} + \B{P}^{(k)}\right) = 0. 
\end{equation}
Here $k$ denotes the index of the Newton's iteration and the solution to \eqref{eq:Newton} is approximated using the right pseudo-inverse of $\D{G}$
\begin{equation*}
\B{u}^{(k+1)} = \B{u}^{(k)} + \B{P}^{(k)}, \quad \B{P}^{(k)}=-\D{G}(\B{u}^{(k)})^\dag G(\B{u}^{(k)}) = -\D{G}^T (\D{G} \, \D{G}^T)^{-1}G 
\end{equation*}
with $\B{u}^{(0)}=\B{\X} + \B{\xi}$. Without loss of generality, we can assume that observation operator $H$ is the identity matrix for a proxy of complete observations. 
The function $G(\B{u})$ has a zero for every orbit of the model. 
The Jacobian of $G$ has an $m(N-1)\times mN$ block structure:
\begin{equation}\label{eq:DG}
	\D{G}(\B{u}) = \begin{bmatrix} 
	- \D{F}_0(u_0) & I \\
	& -\D{F}_1(u_1) & I \\
	& & \ddots & \ddots \\
	& & & -\D{F}_{N-1}(u_{N-1}) & I \end{bmatrix}.
\end{equation}
The Jacobian appears only when acting on a given vector (unit vector for example), and therefore it could be efficiently approximated by finite differences. 
Thus we use an approximation $\D{F} (u)v \approx 1/\varepsilon(F(u+\varepsilon v)-F(u))$.

\section{Shadowing-based DA method}\label{sec:SDA}
In this section, we assume that the observation operator $H$ is not the identity matrix. Therefore, we assume that an initial guess for a shadowing-type DA method is
\begin{equation}\label{eq:IG}
	\B{u}^{(0)} = H^T \B{y}+ (I-H^T H)\B{x}^{\rm b}, 
\end{equation}
where $\B{x}^{\rm b}$ is a so-called background trajectory---a solution of~\eqref{eq:map} with an arbitrary initial condition.

We seek an update $\B{\Pi}^{(k)}$ by approximately solving
\begin{equation*}\label{eq:CG}
	G \left(\B{u}^{(k)} + \B{\Pi}^{(k)}\right) = 0. 
\end{equation*}
using the Levenberg-Marquardt regularization
\begin{equation}\label{eq:regLM}
\B{u}^{(k+1)}=\B{u}^{(k)}+\B{\Pi}^{(k)}, \quad \B{\Pi}^{(k)}=-\Sigma\D{G}^T \left(\D{G} \Sigma \D{G}^T + \alpha^{(k)} Q\right)^{-1}G,
\end{equation}
where $G = G(\B{u}^{(k)})$ defined in~\eqref{eq:residual}, $\D{G} = \D{G}(\B{u}^{(k)})$ defined in~\eqref{eq:DG}, $Q$ is a given positive definite matrix,
and $\alpha^{(k)}>0$.
Here 
\begin{equation}\label{eq:Sigma}
 \Sigma:=H^TRH + (I -H^TH)W(I -H^TH),
\end{equation}
where $W$ is a positive definite matrix that has an $mN\times mN$ block diagonal structure $W = {\rm blockdiag}(W_1,\dots,W_N)$.
The matrix $R$ has a $dN\times dN$ block diagonal structure $R = {\rm blockdiag}(R_1,\dots,R_N)$, 
$Q$ has an $m(N-1)\times m(N-1)$ block diagonal structure $Q = {\rm blockdiag}(Q_1,\dots,Q_N)$,
and $H$ has a $dN\times mN$ block diagonal structure $H = {\rm blockdiag}(H_1,\dots,H_N)$,

To bring a parallel to 4DVar, the solution $\B{\Pi}^{(k)}$ to \eqref{eq:regLM} is a minimizer of a cost function
\begin{equation*}
\frac{1}{2}\left[G\left(\B{u}^{(k+1)}\right)\right]^TQ^{-1}\left[G\left(\B{u}^{(k+1)}\right)\right]+\frac{\alpha^{(k)}}{2} \left[\B{\Pi}^{(k)}\right]^T \Sigma^{-1}\left[\B{\Pi}^{(k)}\right].
\end{equation*}

\subsection{Local convergence}
We define a manifold $\mathcal{M}$ by $ \mathcal{M} = \{\B{u}:\ G(\B{u}) = 0$\} and define $\phi$ as
\[
\phi = \B{u}-\Sigma\D{G}^T \left(\D{G} \Sigma \D{G}^T + \alpha Q\right)^{-1}G.
\]
We note that
\begin{equation}\label{eq:Dphi}
 D\phi = I -\Sigma\D{G}^T \left(\D{G} \Sigma \D{G}^T + \alpha Q\right)^{-1}\D{G}\quad\mbox{for}\quad \B{u}\in\mathcal{M}.
\end{equation}
Since $\mathcal{M}$ is a manifold, we define tangent and normal space of $\mathcal{M}$ at $u$ as $\mathcal{T}_{u}\mathcal{M}$ and $\mathcal{N}_{u}\mathcal{M}$, respectively. 
We have $\mathcal{T}_{u}\mathcal{M}\perp \mathcal{N}_{u}\mathcal{M}$ and 
$\mathcal{T}_{u}\mathcal{M}={\rm kern}\left(\Sigma\D{G}^T \left(\D{G} \Sigma \D{G}^T + \alpha Q\right)^{-1}\D{G} \right)$ for $\B{u}\in\mathcal{M}$. 

\begin{lemma}\label{lemma:uniq}
 $\mathcal{M}$ is a set of fixed points for $\phi$ and there is no further fixed points in the vicinity of $\mathcal{M}$. 
\end{lemma}

\begin{theorem}\label{theorem:conv}
 Suppose $\mathcal{M}$ is compact and contained in an open set $\mathcal{U}$. 
 Furthermore, suppose $D\phi$ is continuous in $\mathcal{U}$ and $\|D\phi\rvert_{\mathcal{N}_{u}\mathcal{M}}\|<1$ for all $\B{u}\in\mathcal{M}$. 
 Then the sequence $\B{u}^{(k)} = \phi^k (\B{u}^{(0)})$ converges for $k\to \infty$ to a point on $\mathcal{M}$ 
 if $\B{u}^{(0)}$ is sufficiently near to $\mathcal{M}$. 
\end{theorem}
For proof of both Lemma~\ref{lemma:uniq} and Theorem~\ref{theorem:conv} we refer to~\cite{BrPar01},
where local convergence for a class of general iterative schemes was proven. 

Now we can prove a local convergence result for the shadowing-based DA method~\eqref{eq:regLM}.
First, we define $\Omega=\D{G}^T Q^{-1}\D{G}$.
\begin{lemma}\label{lemma:Shconv}
Suppose $\Sigma$ and $\Omega$ commute. Furthermore, suppose a positive $\alpha$ satisfies $\alpha>\lambda_{\max} 
(\Sigma \Omega\rvert_{\mathcal{N}_{u}\mathcal{M}})/2 - \lambda_{\min} (\Sigma \Omega\rvert_{\mathcal{N}_{u}\mathcal{M}})$. Then
$\|D\phi\rvert_{\mathcal{N}_{u}\mathcal{M}}\|<1$ for all $\B{u}\in\mathcal{M}$. 
\end{lemma}
\begin{proof}
Using the Sherman-Morrison-Woodbury matrix inversion formula~\cite{GoLo96} and assuming that $\alpha\neq 0$, we can rewrite \eqref{eq:Dphi} as
\begin{equation}\label{eq:DphiW}
 D\phi = I - \Sigma \Omega[ \alpha I+ \Sigma \Omega]^{-1},
\end{equation}
where we drop the iteration notation.
Since $\Sigma$ and $\Omega$ commute, $D\phi$ is symmetric. For symmetric matrices norm is equal to spectral radius. 
Thus $\|D\phi\| = \lambda_{\max}(D\phi)$, where $\lambda_{\max}$ denotes maximum eigenvalue.

A maximum eigenvalue of $D\phi$ is 
\[
 \lambda_{\max}(D\phi) = \max\{|1-\lambda_{\max}(\Sigma \Omega [\alpha I+\Sigma \Omega]^{-1})|,|1-\lambda_{\min}(\Sigma \Omega [\alpha I+\Sigma \Omega]^{-1})|\}.
\]
Moreover, 
\[
0\leq\lambda_{\max}(\Sigma \Omega [\alpha I+\Sigma \Omega]^{-1})\leq\lambda_{\max}(\Sigma \Omega)\lambda_{\max}([\alpha I+\Sigma \Omega]^{-1})=
\frac{\lambda_{\max}(\Sigma \Omega)}{\lambda_{\min}(\alpha I+\Sigma \Omega)} = \frac{\lambda_{\max}(\Sigma \Omega)}{\alpha +\lambda_{\min}(\Sigma \Omega)}.
\]
By choosing $\alpha$ such that
\[
 \frac{\lambda_{\max}(\Sigma \Omega)}{\alpha +\lambda_{\min}(\Sigma \Omega)}<2,
\]
we have $|1-\lambda_{\max}(\Sigma \Omega [\alpha I+\Sigma \Omega]^{-1})|<1$ for $\lambda_{\max}(\Sigma \Omega [\alpha I+\Sigma \Omega]^{-1})>0$.

Furthermore, 
\[
0\leq\lambda_{\min}(\Sigma \Omega [\alpha I+\Sigma \Omega]^{-1})\leq\lambda_{\max}(\Sigma \Omega [\alpha I+\Sigma \Omega]^{-1})<2.
\]
Thus we have $|1-\lambda_{\min}(\Sigma \Omega [\alpha I+\Sigma \Omega]^{-1})|<1$ for $\lambda_{\min}(\Sigma \Omega [\alpha I+\Sigma \Omega]^{-1})>0$.

From~\eqref{eq:DphiW} it follows that $\mathcal{T}_{u}\mathcal{M}={\rm kern}\left(\Sigma \Omega [\alpha I + \Sigma \Omega]^{-1}\right)$ for 
$\B{u}\in\mathcal{M}$. Since $\mathcal{T}_{u}\mathcal{M}\perp \mathcal{N}_{u}\mathcal{M}$, we have
\[
 \lambda\left(\Sigma \Omega [\alpha I+\Sigma \Omega]^{-1}\rvert_{\mathcal{N}_{u}\mathcal{M}}\right)>0.
\]

Therefore by choosing $\alpha>\lambda_{\max} (\Sigma \Omega\rvert_{\mathcal{N}_{u}\mathcal{M}})/2 - \lambda_{\min} (\Sigma \Omega\rvert_{\mathcal{N}_{u}\mathcal{M}})$,
we have $\|D\phi\rvert_{\mathcal{N}_{u}\mathcal{M}}\|< 1$ for all $\B{u}\in\mathcal{M}$.
\end{proof}

\begin{corollary}
The sequence $\B{u}^{(k)} = \phi^k (\B{u}^{(0)})$ defined in~\eqref{eq:regLM} converges for $k\to \infty$ to a point on $\mathcal{M}$ 
if $\B{u}^{(0)}$ is sufficiently near to $\mathcal{M}$. 
\end{corollary}
\begin{proof}
 The proof directly follows from Theorem~\ref{theorem:conv} and Lemma~\ref{lemma:Shconv}.
\end{proof}

\begin{corollary}
Suppose $G$ has only one zero. Then for the sequence defined in~\eqref{eq:regLM}
and a final iteration $K$, $\B{u}^{(K)} = \B{\X}$.
\end{corollary}
This rather trivial corollary shows that the shadowing-based DA method converges to the true solution for linear models or convex $G^T G$.
Existence of several zeros of $G$ is equivalent to the problem of several minima of $G^T G$.

We proved local convergence of the algorithm to the solution manifold. We are unable to provide any results on error bounds with respect to the true solution. However, we provide a necessary condition for an estimate to remain in the trust region of observations. This result is useful since a background trajectory $\B{x}^{\rm b}$ has larger error with respect to the truth than observations $\B{y}$. We recall that an initial guess~\eqref{eq:IG} for the algorithm consists of $\B{x}^{\rm b}$ and $\B{y}$.
Then for a good estimate of the true solution, while updating unobserved variables, observed variables need to have the Gauss-Newton updates that are \emph{inside} the trust region of the observations $\B{y}$.

Before we state the result, let us rewrite the shadowing-based DA method in the limit of continuous algorithmic time step. 
Assume we can set $\alpha^{(k)}=\alpha^{(0)}$ for all $k$. Then we introduce notation $h = 1/\alpha$ and rewrite \eqref{eq:regLM} in terms of $h$ 
\begin{equation*}
\B{u}^{(k+1)}=\B{u}^{(k)}+h\B{\Psi}^{(k)}, \quad \B{\Psi}^{(k)}=-\Sigma\D{G}^T \left(h\D{G} \Sigma \D{G}^T +  Q\right)^{-1}G.
\end{equation*}
Then taking the limit of $h\to 0$, we get on $\tau\in[0\ 1]$
\begin{equation}\label{eq:shadensode}
\frac{d\B{u}}{d\tau}=\psi(\B{u}),\quad\mbox{with}\quad \psi(\B{u})=-\Sigma\D{G}^T(\B{u})Q^{-1}G(\B{u}),\quad\mbox{and}\quad \B{u}(0)=\B{u}^0.
\end{equation}
Defining $\Phi(\B{u})=\|Q^{-1/2}G(\B{u})\|^2/2$, the ODE~\eqref{eq:shadensode} becomes
\begin{equation}\label{eq:gradflow}
\frac{d\B{u}}{d\tau}=-\Sigma\nabla\Phi(\B{u}).
\end{equation}
This is a preconditioned gradient descent for $\Phi(\cdot)$ with a preconditioner $\Sigma$. 
We recall that $\Sigma$ is composed of observation covariance matrix $R$ and weighting matrix $W$~\eqref{eq:Sigma}.
We define $H^{\perp} = (I - H^TH)$. 

\begin{lemma}\label{lemma:trust}
Suppose $\|\nabla\Phi(\B{u})\|<1$. Furthermore, suppose $\| H^{\perp}\B{u}(1) - H^{\perp}\B{x}^{\rm b}\|_W< \varepsilon$ for a small positive $\varepsilon$. Then $\| H\B{u}(1) - \B{y}\|_{R}<1-\varepsilon$. 
\end{lemma}

\begin{proof}
By multiplying \eqref{eq:gradflow} with either $H$ or $H^{\perp}$, taking integral from 0 to 1, and then taking the L2-norm, we have
\begin{equation*}
\| H\B{u}(1) - \B{y}\| = \| R \int_{ 0}^1 H\nabla\Phi(\B{u}) d\tau\|, \quad \mbox{and} \quad
\| H^{\perp}\B{u}(1) - H^{\perp}\B{x}^{\rm b}\| = \| W \int_{ 0}^1 H^{\perp}\nabla\Phi(\B{u}) d\tau\|.
\end{equation*}
Due to assumption $\| H^{\perp}\B{u}(1) - H^{\perp}\B{x}^{\rm b}\|_W< \varepsilon$, we have $\| \int_{ 0}^1H^{\perp}\nabla\Phi(\B{u}) d\tau\|<\varepsilon$. This implies that $\| \int_{ 0}^1H\nabla\Phi(\B{u}) d\tau\|<1-\varepsilon$ for a convergent algorithm $\|\nabla\Phi(\B{u})\|<1$ since $H\perp H^\perp$.
In turn, inequality $\| \int_{ 0}^1H\nabla\Phi(\B{u}) d\tau\|<1-\varepsilon$ implies  $\| H\B{u}(1) - \B{y}\|_{R}<1-\varepsilon$,
 and the estimate  $H\B{u}$ consequently remains in the trust region of observations $\B{y}$. 
\end{proof}

\subsection{Existing shadowing-type DA methods}\label{sec:DAmethods}
Now we point out differences between the shadowing-based DA method introduced in this paper and the existing shadowing-type DA methods of \cite{BrPar01,DuSm14}, and of \cite{Leetal18}. We write down the methods in terms of function $\phi$:
\begin{eqnarray*}
\phi^{\cite{BrPar01}}: &= & u - \D G^T \Lambda^{-1} G, \\
\phi^{\cite{DuSm14}}: & = & u - \gamma \D G^T  G, \\
\phi^{\cite{Leetal18}}: &= & u - \D G^T (\D G \D G^T)^{-1} G, \\
\phi^{\rm{rSh}}: &= & u - \Sigma\D G^T (\D G \Sigma\D G^T + \alpha Q)^{-1} G. 
\end{eqnarray*}
In $\phi^{\cite{BrPar01}}$, $\Lambda$ is chosen to be the Laplace operator. It is stated that the choice of $\Lambda$
has great influence on the convergence, though without a rigorous answer whether the Laplace operator is a good choice. Local convergence is proven for the method as for a class of general iterative schemes. 
In $\phi^{\cite{DuSm14}}$, an algorithmic time step $\gamma$ is chosen by tuning. For sufficiently small $\gamma$ 
 convergence of the damped Gauss-Newton method is guaranteed but the convergence rate might be linear~\cite{GoLo96}. In $\phi^{\cite{Leetal18}}$, the convergence rate is quadratic due to the Gauss-Newton method but the local nature of the Gauss-Newton method requires a good initial guess for convergence---thus (a proxy of) completed observations. In $\phi^{\rm{rSh}}$, on the one hand lower bound on $\alpha$ guarantees local convergence but on the other hand the preconditioner $\Sigma$ might deteriorate the convergence rate. The preconditioner $\Sigma$, namely $W$, is required for a good estimation of the true solution.\par

\section{Numerical experiments\label{sec:numerics}}
We note that if $\Sigma = \epsilon I$, then $\Sigma$ and $\Omega$ commute. For partially-observed models, however, 
 $\Sigma$ and $\Omega$ might not commute. Moreover, in practice $\mathcal{N}_u\mathcal{M}$ is not available. 
Therefore, we assume $\lambda_{\rm max}(\Delta t^2\Sigma\Omega) > \lambda_{\rm max}(\Sigma\Omega|_{\mathcal{N}_u\mathcal{M}})$, where $\Delta t$ is time step of a numerical discretization. Furthermore, we assume 
$ \lambda_{\rm min}(\Sigma\Omega|_{\mathcal{N}_u\mathcal{M}})>0$. The latter assumption is fulfilled if projection onto $\mathcal{N}_u\mathcal{M}$ is defined in terms of $\D{G}$.
Then according to Lemma~\ref{lemma:Shconv} we can choose $\alpha$
 \begin{equation}\label{eq:alpha}
\alpha=\Delta t^2\lambda_{\max} 
(\Sigma \Omega)/2.
\end{equation}
Numerical experiments show that choosing such an $\alpha$ 
provides convergence to the manifold $\mathcal{M}$. However, we do not have a rigorous answer whether the assumption $\lambda_{\rm max}(\Delta t^2\Sigma\Omega) > \lambda_{\rm max}(\Sigma\Omega|_{\mathcal{N}_u\mathcal{M}})$ is fulfilled.
 
When computing $\alpha$, we split the eigenvalue problem over one window length $N$ in $N$ eigenvalue problems over $N$ windows length 1. Then in~\eqref{eq:alpha} we use maximum eigenvalue over $N$ windows.
Moreover, to save computational costs we compute  $\alpha$ for an initial guess $\B{u}^{(0)}$ only and fix the same   $\alpha$ throughout the iteration. The maximum number of iteration is 100.
 Model error is chosen to be $Q = 10^{-3}I$. Other values such as $10^{-2}$ and $10^{-4}$ provide equivalent results to $10^{-3}$. We define the weighting matrix $W = w^2I$ in the preconditioner $\Sigma$ and perform sensitivity analysis in terms of $w$. \par

We compare the shadowing-based DA method to WC4DVar and PDA. PDA is initialised at an initial guess $\B{u}^{(0)}$ and an algorithmic time step is chosen as in~\cite{DuSm14}, namely $\gamma=0.1$. The maximum number of iterations is 100. We note that in~\cite{DuSm14} the maximum number of iterations is 1024. However, we keep the same number of iterations 100 for all DA methods. 

Both the shadowing-based DA method and PDA provide an estimation at observation times only. Therefore we use an estimation at observation times as initial condition for forward model propagation to have an estimation at every time step of numerical discretization. 
 
WC4DVar is initialised at a background trajectory $\B{x}^{\rm b}$. The minimization of a cost function is done by a Matlab built-in Levenberg-Marquardt algorithm and stopping when the relative change in the 
cost function compared to the initial value is less then $10^{-6}$ unless 100 iterations are reached. Model error for WC4DVar is $10^{-2}I$, and the background covariance matrix is the identity.

In order to check robustness of the results,
we perform 100 numerical experiments with different realizations of truth $\B{\X}$, observations $\B{y}$, and background trajectory $\B{x}^{\rm b}$.

To analyze the shadowing-based DA method and compare it to other methods, we compute mean over time of $G$-error
\begin{equation}\label{eq:errG}
\mathbb{E}^G = \frac{1}{N}\sum_{n=0}^{N-1}G_n^T G_n,
\end{equation}
mean over time of estimation error with respect to the truth of observed variables
\begin{equation}\label{eq:MerrO}
\mathbb{E}^{\rm O}= \frac{1}{N}\sum_{n=0}^{N-1}\mathbb{E}^{\rm O}_n,
\end{equation}
and of non-observed variables
\begin{equation}\label{eq:MerrNo}
\mathbb{E}^{\rm N}= \frac{1}{N}\sum_{n=0}^{N-1}\mathbb{E}^{\rm N}_n.
\end{equation}
Here errors $\mathbb{E}^{\rm O}_n$ and $\mathbb{E}^{\rm N}_n$ defined as
\begin{equation}\label{eq:errO}
\mathbb{E}^{\rm O}_n= \frac{1}{{\rm rank}(H)}(H u_n - H\X_n)^T (H u_n - H \X_n)
\end{equation}
and 
\begin{equation}\label{eq:errNo}
\mathbb{E}^{\rm N}_n= \frac{1}{{\rm rank}(I-H^TH)}[(I-H^TH)(u_n - \X_n)]^T [(I-H^TH)(u_n - \X_n)],
\end{equation}
respectively, and $n$ is an index of numerical time step not observation time step.
We also compute a cost function with respect to observations
\begin{equation}\label{eq:cost}
C= \frac{1}{(k_2-k_1+1) {\rm rank}(H)}\sum_{k=k_1}^{k_2}(H u_k - y_k)^T (H u_k - y_k),
\end{equation}
where $k$ is an index of observation time step, $k_1\geq 0$, and  $k_2\leq N-1$.

\subsection{Application to the Lorenz 63 model}
The well-known Lorenz attractor \cite{Lo63} is a chaotic dynamical system commonly used as a test problem for data assimilation algorithms. The L63 model is
\begin{equation}\label{eq:L63ode}
    \dot{x}^1 = \sigma(x^2-x^1), \quad 
    \dot{x}^2 = x^1(\rho-x^3)-x^2, \quad 
    \dot{x}^3 = x^1x^2-\beta x^3,
\end{equation}
where $\sigma=10$, $\beta=\frac{8}{3}$ and $\rho=28$.  
The differential equations are discretized with a forward Euler scheme with time step $\Delta t=0.005$. (We have also considered Runge-Kutta 4th order but since it gives similar results, it is omitted in the paper.) 
We generate a set of observations computing a trajectory of L63 on $t\in[0,100]$, with a spin-up of $[-25,0]$ for a true trajectory to reside on the attractor.
Observations are obtained by perturbing a reference (true) trajectory with random Gaussian iid noise with zero mean and covariance $R=8I$. 
The observations of $x^1$-variable only are drawn every $\Delta t_{\rm obs}=0.05$.
Then the map $F_n$~\eqref{eq:map} corresponds to 10 forward Euler steps.
This map is used to define $G$ and the derivatives of this map are needed for the shadowing iteration. 
The assimilation windows is $\Delta t_{\rm ass}=5$. \par

\begin{figure}[ht]
	\centering
{\includegraphics[width=\textwidth]{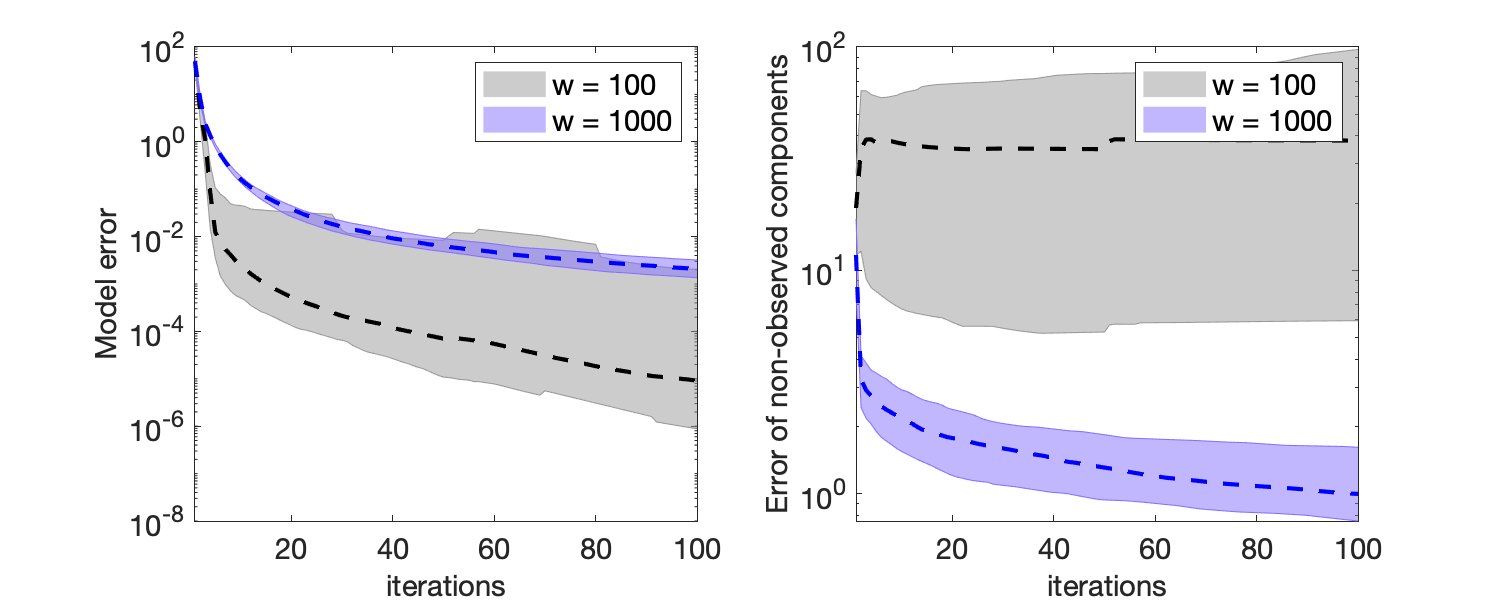}}
	\caption{Application to L63. Error of the shadowing-based DA method as a function of iterations: median (dashed line) +/- one standard deviation (shadowed area) over 100 simulations. In grey error is shown for weighting matrix $w=100$, in blue for $w=1000$.
		On the left: mean over time of $G$-error. On the right: mean over time of error with respect to the truth of non-observed variables.}\label{fig:Fig63_1}
\end{figure}

\begin{figure}[ht]
	\centering
{\includegraphics[width=\textwidth]{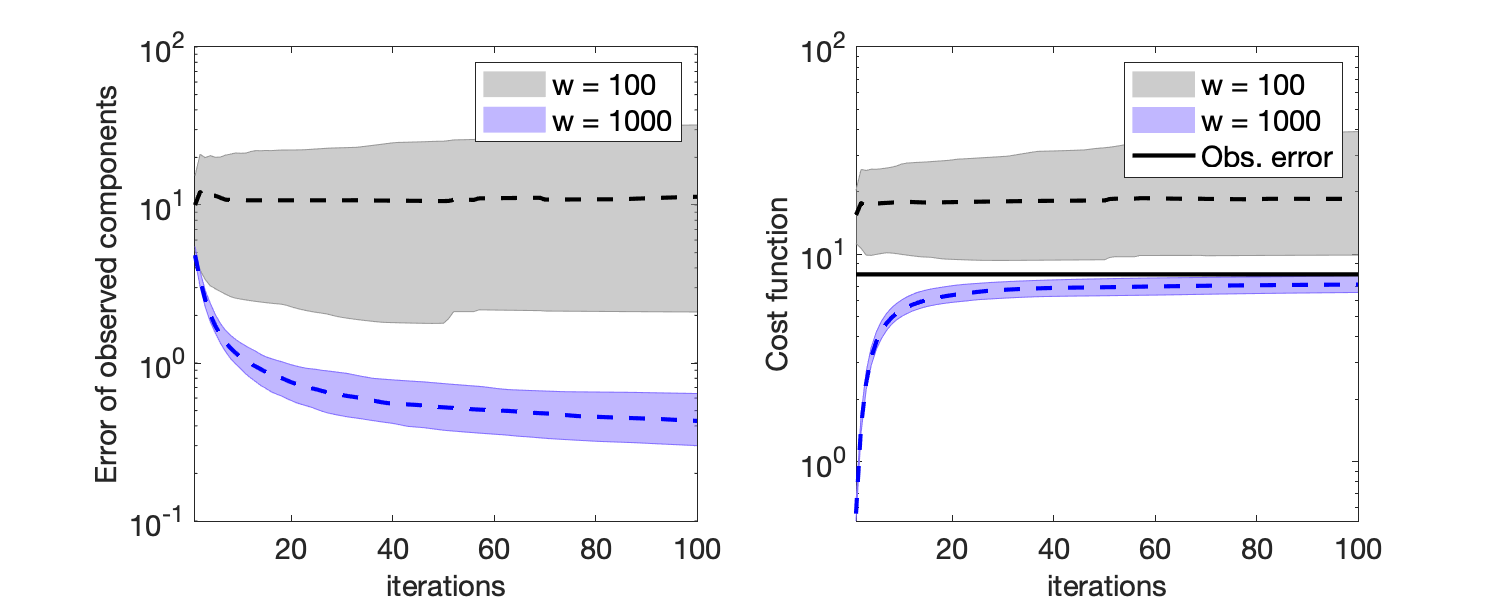}}
	\caption{Application to L63. Error of the shadowing-based DA method as a function of iterations: median (dashed line) +/- one standard deviation (shadowed area) over 100 simulations. In grey error is shown for weighting matrix $w=100$, in blue for $w=1000$.
On the left: mean over time of error with respect to the truth of observed variables. On the right: mean over time of cost function of observed variables.}\label{fig:Fig63_2}
\end{figure}

In Figure~\ref{fig:Fig63_1} we display $G$-error~\eqref{eq:errG} on the left and error with respect to the truth of non-observed variables~\eqref{eq:MerrNo} on the right as a function of iteration. We remark that small $w=100$ gives quicker convergence to the manifold $\mathcal{M}$, while large $w=1000$ requires more iterations to reach the same error on average. However, error with respect to the truth of non-observed variables is decreasing over iteration for large $w=1000$, while increasing for small $w=100$.
In Figure~\ref{fig:Fig63_2} we plot error with respect to the truth of observed variables~\eqref{eq:MerrO} on the left 
and cost function~\eqref{eq:cost} on the right as a function of iteration, where solid black line is for observation error. We see again that large $w=1000$ gives better estimation of observed variables than small $w=100$. \par

When analyzing the cost function, we see that for small $w=100$ the cost function quickly underestimates the observation error. In inverse problems this phenomenon is often referred as observations overfitting, though a cost function is there decreasing not increasing and the observation error is overestimated not underestimated, see e.g.~\cite{Ha97}. For the shadowing-based DA method the cost function~\eqref{eq:cost} at the first iteration is zero, because the algorithm is initialized at $\B{u}^{(0)}$~\eqref{eq:IG}. The cost function increases over iteration due to a search for a noise-free orbit. When the cost function is larger than the observation error $R$, an estimate is not in a ball of radius $R$ centred at the true trajectory, resulting in a larger error with respect to the truth. Therefore, we need to prevent the cost function becoming larger than $R$. A classical approach in inverse problems is to stop the iteration when this occurs. In the shadowing-based DA method this approach is questionable due to cost function increasing over iteration. Instead, we propose to tune the preconditioner $\Sigma$~\eqref{eq:Sigma}, namely the weighting matrix $W$, to obtain the correct behaviour of the cost function. We see that the large value of $w=1000$ results in the cost function approaching the observation error from below. This is an indication of correctly tuned $w$. Thus the role of preconditioner $\Sigma$ is to keep descend steps in the direction of observed variables small compared to descend steps in the direction of non-observed variables.
As the iteration proceeds, observed variables get denoized as well and the algorithm finds a (pseudo-)orbit compatible with observations. We would like to stress that the cost function~\eqref{eq:cost} depends only on observations, not the truth.\par

\begin{figure}[ht]
	\centering
{\includegraphics[width=\textwidth]{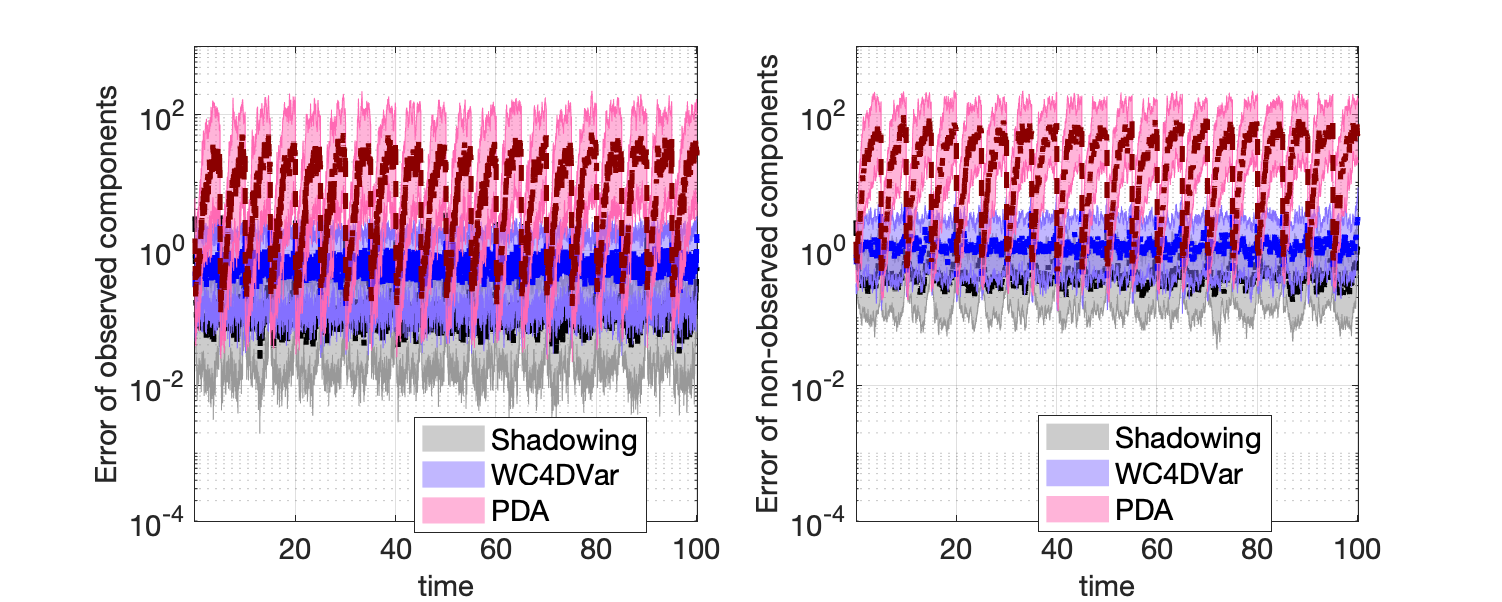}}
	\caption{Application to L63. Error as a function of time: median (dashed line) +/- one standard deviation (shadowed area) over 100 simulations. On the left: error with respect to the truth of observed variables. On the right: error with respect to the truth of non-observed variables. The shadowing-based DA method with $w=1000$ in grey, WC4DVar in blue, and PDA in pink.}\label{fig:Fig63_3}
\end{figure}

In Figure~\ref{fig:Fig63_3} we compare the shadowing-based DA method with $w=1000$ to WC4DVar and PDA, where we plot error with respect to the truth over time of observed variables~\eqref{eq:errO} and of non-observed variables~\eqref{eq:errNo} on the left and right, respectively. We see that the correct choice of the preconditioner is essential for shadowing-type DA methods, since for fully observed L63 PDA and the shadowing-based DA method perform comparably (not shown) but for partially observed L63 PDA perform poorly. It is also remarkable that the shadowing-based DA method with tuned $w$ outperforms WC4DVar. 

\subsection{Application to the Lorenz 96 model}
Lorenz \cite{Lo96} proposed the following model as an example of a simple one-dimensional model with features of the atmosphere.  The L96 model is
\begin{equation}\label{eq:L96ode}
\dot{x^l}=-x^{l-2}x^{l-1}+x^{l-1}x^{l+1}-x^l+\mathcal{F},\qquad (l=1,...,d),
\end{equation}
where the dimension $d$ and forcing $\mathcal{F}$ are parameters. Cyclic boundary conditions are imposed.
We implement the L96 model with the standard parameter choices $d=36$ and $\mathcal{F}=8$. 
The differential equations are discretized with a forward Euler scheme with time step $\Delta t=0.005$. (We have also considered Runge-Kutta 4th order but since it gives similar results, it is omitted in the paper.) 
We generate a set of observations computing a trajectory of L96 on $t\in[0,100]$, with a spin-up of $[-25,0]$ for a true trajectory to reside on the attractor.
Observations are obtained by perturbing a reference (true) trajectory with random Gaussian iid noise with zero mean and covariance $R=8I$. 
The observations of every 2nd variable are drawn every $\Delta t_{\rm obs}=0.05$.
Then the map $F_n$~\eqref{eq:map} corresponds to 10 forward Euler steps.
This map is used to define $G$ and the derivatives of this map are needed for the shadowing iteration. 
The assimilation windows is $\Delta t_{\rm ass}=5$.

\begin{figure}[ht]
	\centering
{\includegraphics[width=\textwidth]{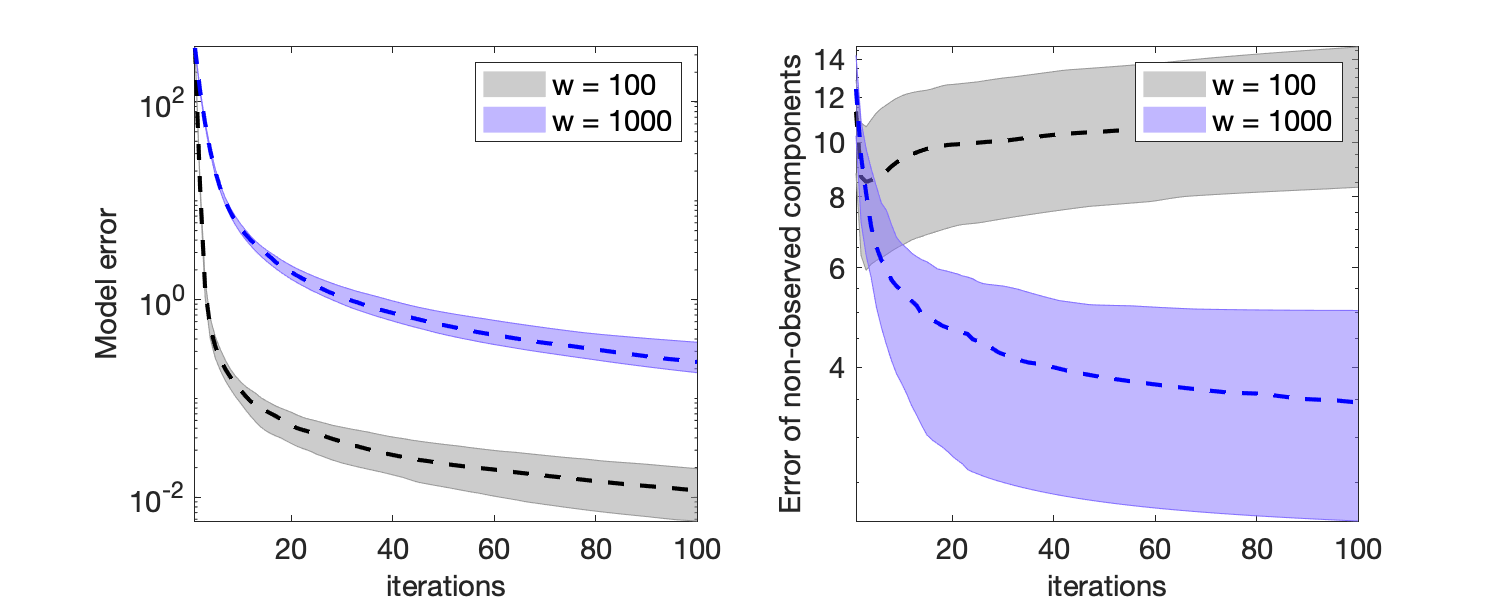}}
	\caption{Application to L96. Error of the shadowing-based DA method as a function of iterations: median (dashed line) +/- one standard deviation (shadowed area) over 100 simulations. In grey error is shown for weighting matrix $w=100$, in blue for $w=1000$.
		On the left: mean over time of $G$-error. On the right: mean over time of error with respect to the truth of non-observed variables.}\label{fig:Fig96_1}
\end{figure}

In Figure~\ref{fig:Fig96_1} we display $G$-error~\eqref{eq:errG} on the left and error with respect to the truth of non-observed variables~\eqref{eq:MerrNo} on the right as a function of iteration. As for L63 displayed in Figure~\ref{fig:Fig63_1}, large $w=1000$ requires more 
iterations to reach the same $G$-error than small $w=100$. Error with respect to the truth of non-observed variables decreases over iteration for large $w=1000$ while increases for small $w=100$.

\begin{figure}[ht]
	\centering
{\includegraphics[width=\textwidth]{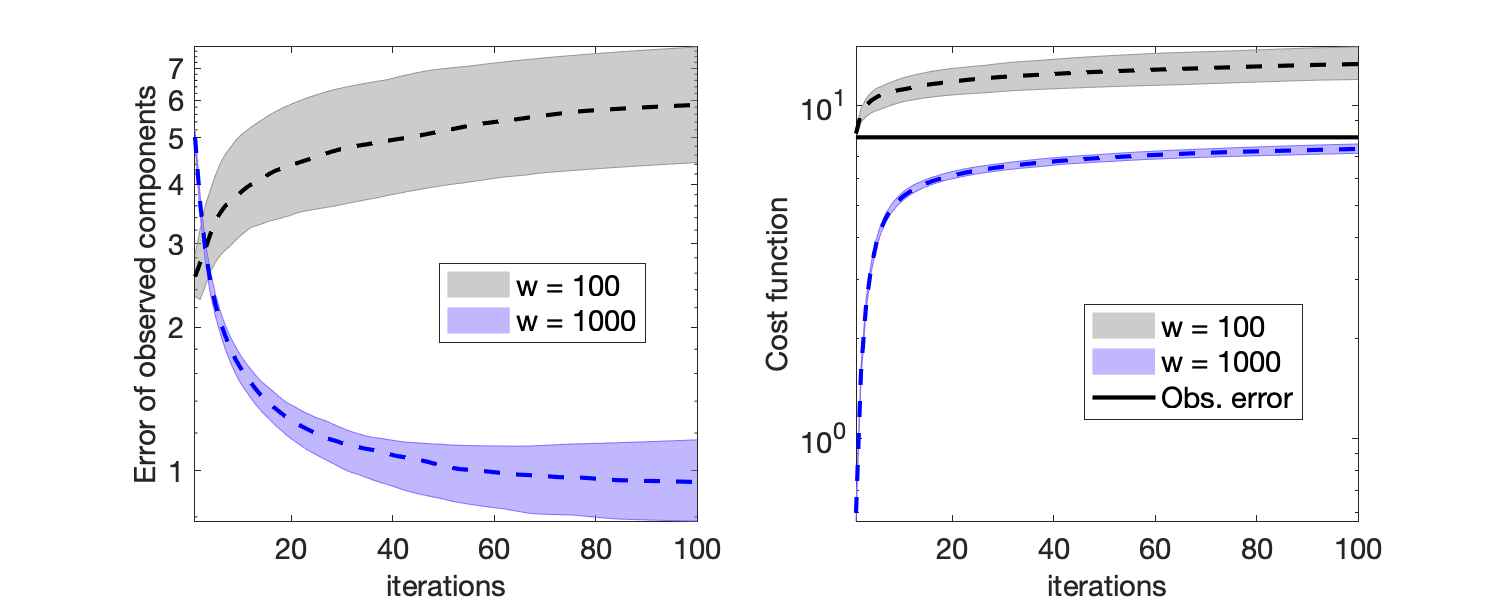}}
	\caption{Application to L96. Error of the shadowing-based DA method as a function of iterations: median (dashed line) +/- one standard deviation (shadowed area) over 100 simulations. In grey error is shown for weighting matrix $w=100$, in blue for $w=1000$.
On the left: mean over time of error with respect to the truth of observed variables. On the right: mean over time of cost function of observed variables.}\label{fig:Fig96_2}
\end{figure}

In Figure~\ref{fig:Fig96_2}, we plot error with respect to the truth of observed variables~\eqref{eq:MerrO} 
and cost function~\eqref{eq:cost} as a function of iteration on the left and on the right, respectively.
A better estimation of observed variables is obtained with large $w=1000$ than with small $w=100$, as was the case for L63 displayed in Figure~\ref{fig:Fig63_2}. Moreover, small $w=100$ gives a considerable increase in the error. The cost function is underestimated with small $w=100$ and well estimated with large $w=1000$. Thus the preconditioner $\Sigma$ with $w=1000$ is optimal.

\begin{figure}[ht]
	\centering
{\includegraphics[width=\textwidth]{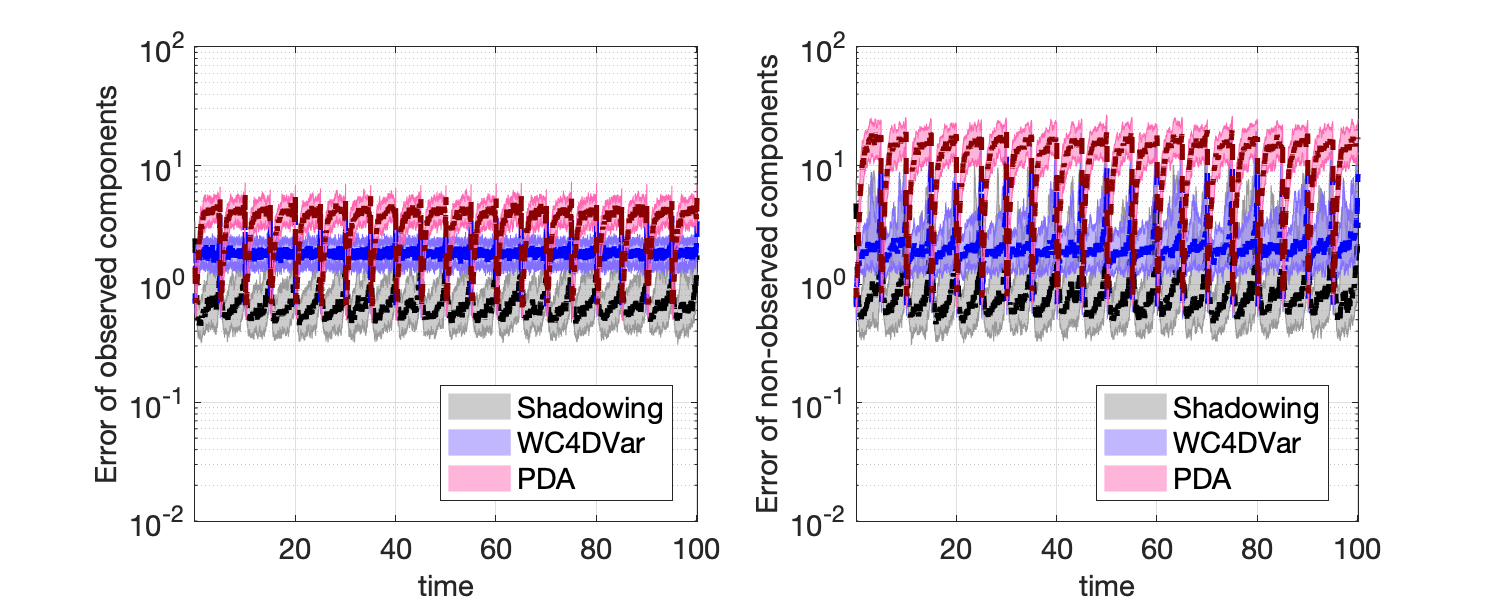}}
	\caption{Application to L96. Error as a function of time: median (dashed line) +/- one standard deviation (shadowed area) over 100 simulations. On the left: error with respect to the truth of observed variables. On the right: error with respect to the truth of non-observed variables. The shadowing-based DA method with $w=1000$ in grey, WC4DVar in blue, and PDA in pink.}\label{fig:Fig96_3}
\end{figure}

In Figure~\ref{fig:Fig96_3} we compare the shadowing-based DA method with $w=1000$ to WC4DVar and PDA, where we plot error with respect to the truth over time of observed variables~\eqref{eq:errO} and of non-observed variables~\eqref{eq:errNo} on the left and right, respectively. 
Here we see that  the shadowing-based DA method with correctly chosen preconditioner $\Sigma$ 
outperforms both WC4DVar and PDA.

\section{Conclusions\label{sec:con}}
We have introduced a shadowing-based DA method for partial observations based on the regularized Gauss-Newton method. We proved local convergence of the method and derived a lower bound for the algorithmic time step required for the method to converge to the manifold $G(u) = 0$. We also introduced a preconditioner for the shadowing-based DA  method. The preconditioner scales the descend steps such that the descend step of non-observed variables is large compared to observed variables. This allows the algorithm to find a solution of $G(u) = 0$ in the vicinity of the truth.
Numerical experiments with the Lorenz 63 and Lorenz 96 models show encouraging results: the shadowing-based DA method outperforms both WC4Var and PDA. The shadowing-based DA method is more expensive than PDA and WC4Var, since it requires finding eigenvalues at the first iteration, forming large matrices and inverting them. Therefore future directions include decreasing computational costs, a rigorous answer to the numerical choice of $\alpha$, and error bounds with respect to the truth.  

\section{Acknowledgements}
This work is part of the research programme Mathematics of Planet Earth 2014 EW with project number 657.014.001, which is financed by the Netherlands Organisation for Scientific Research (NWO).


%
%

%


\bibliographystyle{siamplain}
\bibliography{main}

\end{document}